\documentclass[11pt]{article}
\usepackage{comment}
 \usepackage{tabularx,lscape,longtable}
 \usepackage{comment}
 \usepackage{amsfonts}
 \usepackage{multirow}
 \usepackage{latexsym}
 \usepackage{amssymb}
 \usepackage{amsmath}
 \usepackage{array}
 \usepackage{tikz}
 \usepackage{graphicx}
 \usepackage{tfrupee} 
\begin{document}
\newtheorem{t1}{Theorem}[section]
\newtheorem{d1}{Definition}[section]
\newtheorem{c1}{Corollary}[section]
\newtheorem{l1}{Lemma}[section]
\newtheorem{r1}{Remark}[section]
\newtheorem{e1}{Example}[section]
\title{Errors Due to Departure from Independence in \\Exponential Series System}
\author{Asok K. Nanda\footnote{corresponding author. e-mail: asok.k.nanda@gmail.com}, Sanjib Gayen\\
	Department of Mathematics and Statistics\\Indian Institute of Science Education and Research Kolkata\\West Bengal, India \and Shovan Chowdhury\\
	Quantitative Methods and Operations Management Area\\
	Indian Institute of Management, Kozhikode\\
	Kerala, India.}
\date{\today}
\maketitle
\begin{abstract}
\noindent In reliability and life testing when the exponentially distributed components are put in series, it is generally assumed that the lifetimes of the components are independently distributed, which leads to some errors if they are not actually independent. In this paper, we study the relative errors incurred in different reliability measures due to such assumptions when actually they follow some bivariate exponential distributions.
\end{abstract}
{\bf Key Words \& Phrases:} Hazard rate, mean residual life function, reliability function, reversed hazard rate function.
\section{Introduction}
Consider a two-component series system. For such a system the failure of any one of the two components causes the system to fail. A common assumption in reliability and life testing, in modeling and analyzing data from such a system, is that they are independent and exponentially distributed. Sometimes such an assumption of independence is clearly false. Thus, to make an understanding, a number of bivariate exponential distributions have been obtained by many authors viz. Gumbel ($1960$), Freund ($1961$), Marshall and Olkin ($1967$), Block and Basu ($1974$), Cowan ($1987$) and Sarkar ($1987$) among others. The assumption of independence when in fact the lifetime distribution is a bivariate exponential leads to some error in the analysis of data.\\
\hspace*{.2in} Klein and Moeschberger ($1986$, $1987$) and Moeschberger and Klein ($1984$) have studied the relative errors in some bivariate exponential distributions. Gupta and Gupta ($1990$) have also studied the relative errors in different reliability measures.\\
\hspace*{.2in} In this paper, we study the relative errors in different reliability measures viz. {\it reliability function, failure rate function, mean residual life function} and {\it reversed hazard rate} (RHR) {\it function} (RHR function plays an important role in the analysis of left-censored data) under the assumption of independence when they actually follow either of the bivariate exponential distributions due to Block and Basu ($1974$), Cowan ($1987$), Freund ($1961$), Gumbel ($1960$), Marshall and Olkin ($1967$) and Sarkar ($1987$).\\
\hspace*{.2in} Consider a two-component series system whose components jointly follow bivariate exponential distribution. Let $X_i$ be the lifetime of the $i^{th}$ component ($i=1,2$) of the system. Then the joint distribution of $X_1$ and $X_2$ may be one of the following:
\begin{enumerate}
\item Independent:
$$\bar{F}_1(x_1,x_2)=e^{-\lambda_1x_1-\lambda_2x_2},\quad
\lambda_i,x_i>0,\;i=1,2;$$
\item Gumbel I:
$$\bar{F}_2(x_1,x_2)=e^{-\lambda_1x_1-\lambda_2x_2-\lambda_{12}x_1x_2},\quad
\lambda_{12}\ge 0; ~0 \leq \lambda_{12} \leq \lambda_{1} \lambda_{2}, $$
$\lambda_i,x_i>0,\;i=1,2;$
\item Gumbel II:
$$\bar{F}_3(x_1,x_2)=\left[1+\alpha\left(1-e^{-\lambda_1x_1}\right)\left(1-e^{-\lambda_2x_2}\right)\right]e^{-\lambda_1x_1-\lambda_2x_2},\quad |\alpha|<1,$$
$\lambda_i,x_i>0,\;i=1,2;$
\item Gumbel III:
$$\bar{F}_4(x_1,x_2)=e^{-((\lambda_{1} x_{1})^{m} + ( \lambda_{2} x_{2})^{m} )^{\frac{1}{m}}}, \quad
\lambda_i,x_i>0,\;i=1,2;~m \geq 1 $$
\item Freund:\\
$$\bar{F}_5(x_1,x_2)=\left\{\begin{array}{ll}
\frac{\lambda_1}{\lambda-\theta_2}e^{-(\lambda-\theta_2)x_1-\theta_2x_2}+    \frac{\lambda_2-\theta_2}{\lambda-\theta_2}e^{-\lambda x_2},&x_1\le x_2\\  &\\  \frac{\lambda_2}{\lambda-\theta_1}e^{-(\lambda-\theta_1)x_2-\theta_1x_1}+    \frac{\lambda_1-\theta_1}{\lambda-\theta_1}e^{-\lambda x_1},&x_1> x_2
                                 \end{array}\right.$$
$\lambda_i,\theta_i,x_i>0,\;i=1,2;\lambda=\lambda_1+\lambda_2;$
\item Marshall and Olkin:
$$\bar{F}_6(x_1,x_2)=e^{-\lambda_1x_1-\lambda_2x_2-\lambda_{12}\max(x_1,x_2)},\quad
\lambda_i,x_i>0,\;i=1,2;\;\lambda_{12}> 0;$$
\item Block and Basu:
$$\bar{F}_7(x_1,x_2)=\frac{\lambda^*}{\lambda}e^{-\lambda_1x_1-\lambda_2x_2-\lambda_{12}\max(x_1,x_2)}-\frac{\lambda_{12}}{\lambda}e^{-\lambda^*\max(x_1,x_2)},\quad \lambda_i,x_i>0,\;i=1,2,$$
$\lambda_{12}\ge 0,\lambda=\lambda_1+\lambda_2,\lambda^*=\lambda+\lambda_{12};$
\item Cowan:
$$\bar{F}_8(x_1,x_2)=e^{-\frac{1}{2}\left[\lambda_1x_1+\lambda_2x_2+\left(\lambda_1^2x_1^2+\lambda_2^2x_2^2-2\lambda_1\lambda_2x_1x_2\cos\theta\right)^{1/2}\right]},\quad \theta\in ( 0,\pi],$$
$\lambda_i,x_i>0,i=1,2;$
\item Sarkar:
$$\bar{F}_9(x_1,x_2)=\left\{\begin{array}{ll}                           e^{-\left(\lambda_2+\lambda_{12}\right)x_2}\left[1-\left(1-e^{-\lambda_1x_2}\right)^{-\nu}.\left(1-e^{-\lambda_1x_1}\right)^{1+\nu}\right],&0<x_1 < x_2\\&\\
e^{-\left(\lambda_1+\lambda_{12}\right)x_1}\left[1-\left(1-e^{-\lambda_2x_1}\right)^{-\nu}.\left(1-e^{-\lambda_2x_2}\right)^{1+\nu}\right],&0<x_2\le x_1,
                           \end{array}\right.$$
where $\lambda_i,x_i>0,i=1,2;\lambda_{12}>0,\nu=\lambda_{12}/(\lambda_1+\lambda_2)$.
\end{enumerate}
An error may occur by assuming independence of the components when in fact the distribution of the components is described by one of the above models.\\
\hspace*{.2in} In Section \ref{def} of this paper we give various reliability measures along with the corresponding relative errors. Calculation of errors in different reliability measures corresponding to different bivariate distributions have been studied in Section \ref{cal}, whereas  Section \ref{an} reports the detailed analysis of these errors. It is to be mentioned here that by $a\stackrel{sign}{=}b$ we mean that $a$ and $b$ have the same sign.
\section{Definitions and Preliminaries}\label{def}
Let $T$ be a nonnegative random variable denoting the lifetime of a component having
distribution function $F(\cdot)$, density function $f(\cdot)$, survival function $\bar{F}(\cdot)=1-F(\cdot)$. Then the mean residual life function $e(\cdot)$ is defined as
$$e(t)=E[T-t|T>t]=\int_t^\infty\frac{\bar{F}(x)dx}{\bar{F}(t)},$$
 whereas the failure rate function $r(\cdot)$ and the reversed hazard rate function $\mu(\cdot)$ of $T$ are defined, respectively, as
$r(t)=f(t)/\bar{F}(t)$ and $\mu(t)=f(t)/F(t)$. The relative errors in survival function and in mean residual life function are given respectively as
$$\frac{\bar{F}_D(t)-\bar{F}_I(t)}{\bar{F}_I(t)}\qquad{\rm and}\qquad \frac{e_D(t)-e_I(t)}{e_I(t)},$$
whereas that in failure rate function and in reversed hazard rate function are defined respectively as
$$\frac{r_D(t)-r_I(t)}{r_I(t)}\qquad{\rm and}\qquad \frac{\mu_D(t)-\mu_I(t)}{\mu_I(t)},$$
where $D$ and $I$ stand for dependent and independent models respectively. Throughout
this paper, the words positive (negative) and nonnegative (non-positive) are used interchangeably.

It is well known that the reliability measures described in Table \ref{tab1} are equivalent in the sense that knowing the one other can be uniquely determined from the relationships\\
$$\bar{F}(t)=e^{-\int_0^tr(x)dx},\quad F(t)~=e^{-\int_t^\infty\mu(x)dx},\quad r(t)=\frac{1+e'(t)}{e(t)}.$$
It is to be mentioned here that although the reliability measures are equivalent, the relative errors in these measures do not exhibit similar property as Table~\ref{tab2} shows, and hence separate study for each of the reliability measures is necessary.
\section{Calculation of Different Reliability Measures}\label{cal}
Let the lifetime of a two-component series system be denoted by $T$. Then the survival function of $T$ is given by
\begin{eqnarray*}
\bar{F}_T(t)&=&P(T>t)\\
&=&P(\min\{X_1,X_2\}>t)\\
&=&P(X_1>t,X_2>t)\\
&=&\bar F(t,t),
\end{eqnarray*}
where $\bar F(x,y)=P(X_1>x,X_2>y)$ is the probability that $X_1$ survives for at least $x$ units and $X_2$ survives for at least $y$ units of time.
\subsection{Reliability Measures for Gumbel I Model}
From the expression of the survival function of the Gumbel I distribution we have the survival function, $\bar F_T$, of the two-component series system given by
$$\bar F_T(t)=e^{-\lambda t-\lambda_{12}t^2},$$
where $\lambda=\lambda_1+\lambda_2$. Clearly, the failure rate, $r_T(t)$, of $T$ is given by $\lambda+2\lambda_{12}t$. The mean residual life for the distribution is obtained as
\begin{eqnarray*}
e_T(t)&=&\frac{\int_t^\infty \bar F_T(x)dx}{\bar F_T(t)}\\
&=&\frac{\int_t^\infty e^{-\lambda x-\lambda_{12}x^2}dx}{e^{-\lambda t-\lambda_{12}t^2}}
\end{eqnarray*}
\begin{eqnarray*}
&=&e^{\lambda t + \lambda_{12} t^2} \int_t^\infty e^{- \lambda x - \lambda_{12} x^2} dx \\
&=& e^{\left(\sqrt{\lambda _{12}}t +\frac{\lambda}{2 \sqrt{\lambda_{12} } }\right)^2}	
\int _t^{\infty} e^{-\frac{1}{2}\left(\sqrt{2 \lambda_{12}}x+ \frac{\lambda}{\sqrt{2 \lambda _{12}}}\right)^2} dx \\
&=&e^{\lambda_{12}\left(t+\frac{\lambda}{2\lambda_{12}}\right)^2} \frac{1}{\sqrt{2\lambda_{12}}}\int_{\sqrt{2\lambda_{12}}t+\frac{\lambda}{\sqrt{2\lambda_{12}}}}^{\infty}e^{-\frac{x^2}{2}} dx\\ 
&=&\sqrt{\frac{\pi}{\lambda_{12}}}e^{\delta (t)} \left(\frac{1}{\sqrt{2 \pi}} \int _{- \infty}^{\infty}e^{- \frac{x^2}{2}} dx - \frac{1}{\sqrt{2 \pi}}\int _{- \infty} ^{\sqrt{2 \delta (t)}}e^{- \frac{x^2}{2}}dx\right)\\
&=&\\
&=&\sqrt{\frac{\pi}{\lambda_{12}}}e^{\delta (t)}\left(1- \Phi\left(\sqrt{2 \delta (t)}\right)\right),
\end{eqnarray*}
where $\Phi(\cdot)$ is the distribution function of the standard normal distribution and $\delta(t)=\lambda_{12}\left(t+\frac{\lambda}{2\lambda_{12}}\right)^2$. It is to be noted that the distribution function, $F_T$, of the system is given by 
$$F_T(t)=P(\min\{X_1,X_2\}\leqslant t)=1-e^{-\lambda t-\lambda_{12}t^2}$$
so that the reversed hazard rate function, $\mu_T(\cdot)$, is obtained as
$$\mu_T(t)=\frac{\lambda+2\lambda_{12}t}{e^{\lambda t+\lambda_{12}t^2}-1}.$$
\subsection{Reliability Measures for Gumbel II Model}
From the expression of the survival function of the Gumbel II distribution we have the survival function, $\bar F_T$, of the two-component series system as
$$\bar F_T(t)=e^{-\lambda t}\left[1+\alpha\left(1-e^{-\lambda_1t}\right)\left(1-e^{-\lambda_2t}\right)\right].$$
 After simple calculation we have the failure rate, $r_T(t)$, of $T$ as
$$r_T(t)=\lambda-\frac{h'(t)}{h(t)},$$
where $h(t)=1+\alpha\left(1-e^{-\lambda_1t}\right)\left(1-e^{-\lambda_2t}\right)$. The mean residual life for the distribution is obtained as
\begin{eqnarray*}
e_T(t)&=&\frac{\int_t^\infty \bar F_T(x)dx}{\bar F_T(t)}\\
&=&\frac{\int_t^\infty h(x)e^{-\lambda x}dx}{h(t)e^{-\lambda t}}\\
&=&\frac{e^{\lambda t}}{h(t)}\int_t^\infty \left[(1+\alpha)e^{-\lambda x}+\alpha \left(e^{-2\lambda x}-e^{-\lambda x -\lambda_1x}-e^{-\lambda x -\lambda_2 x}\right)\right]dx \\
&=&\frac{1}{h(t)}\left[\frac{(1+\alpha)}{\lambda}-\alpha \left(\frac{e^{-\lambda_1 t}}{\lambda + \lambda_1}+\frac{e^{-\lambda_2 t}}{\lambda + \lambda_2}-\frac{e^{-\lambda t}}{2 \lambda}\right)\right]\\
&=&\frac{1}{h(t)}\left[\frac{1+\alpha}{\lambda} - \alpha g(t)\right],
\end{eqnarray*}
where $g(t)=\frac{e^{-\lambda_1t}}{\lambda_1+\lambda}+
\frac{e^{-\lambda_2t}}{\lambda+\lambda_2}-\frac{e^{-\lambda t}} {2\lambda}$.
It is to be noted that the density function, $f_T$, of the system is given by 
$$f_T(t)=e^{-\lambda t}\left(\lambda h(t)-h'(t)\right)$$
so that the reversed hazard rate function, $\mu_T(\cdot)$, is obtained as
$$\mu_T(t)=\frac{\lambda h(t)-h'(t)}{e^{\lambda t}-h(t)}$$
\subsection{Reliability Measures for Gumbel III Model}
Note that the survival function, $\bar F_T$, of the two-component series system of Gumbel III model is given by 
$$\bar F_T(t)=e^{-t(\lambda_1^m + \lambda_2^m )^{1/m}},$$
which gives the failure rate, $r_T(t)$, and the mean residual life, $e_T(t)$, of $T$ as
$$r_T(t)= (\lambda_1^m + \lambda_2^m )^{1/m}\qquad {\rm and}\qquad e_T(t)= (\lambda_1^m + \lambda_2^m )^{-1/m}.$$
The reversed hazard rate function, $\mu_T(\cdot)$, is obtained as
$$\mu_T(t)=\frac{(\lambda_1^m + \lambda_2^m )^{1/m}}{e^{t(\lambda_1^m + \lambda_2^m )^{1/m}}-1}.$$
\subsection{Reliability Measures for Freund's Model}
Since the survival function of Freund model is given by $\bar F_T(t)=e^{-t(\lambda_1+\lambda_2)}$, we immediately get the failure rate, $r_T(t)$, and the mean residual life, $e_T(t)$, of $T$ as
$$r_T(t)= \lambda_1+\lambda_2\qquad {\rm and}\qquad e_T(t)= \frac{1}{\lambda_1+\lambda_2}.$$
 The reversed hazard rate function, $\mu_T(\cdot)$, is obtained as
$$\mu_T(t)=\frac{\lambda_1+\lambda_2}{e^{t(\lambda_1+\lambda_2)}-1}$$
\subsection{Reliability Measures for Cowan's Model}
The survival function of Cowan's model is given by 
$$\bar F_T(t)=e^{-\frac{t}{2}\left(\lambda_1+\lambda_2 +\sqrt{\lambda_1^2+\lambda_2^2-2\lambda_1\lambda_2\cos\theta}\right)}.$$
Immediately we get the failure rate, $r_T(t)$, and the mean residual life, $e_T(t)$, of $T$ as
$$r_T(t)= \frac{1}{2}\left[\lambda_1+\lambda_2 +\sqrt{\lambda_1^2+\lambda_2^2-2\lambda_1\lambda_2\cos\theta}\right]$$
 and
$$e_T(t)= 2\left(\lambda_1+\lambda_2 +\sqrt{\lambda_1^2+\lambda_2^2-2\lambda_1\lambda_2\cos\theta}\right)^{-1}.$$
The reversed hazard rate function, $\mu_T(\cdot)$, is obtained as
$$\mu_T(t)=\frac{\lambda_1+\lambda_2 +\sqrt{\lambda_1^2+\lambda_2^2-2\lambda_1\lambda_2\cos\theta}}{2\left(e^{\frac{t}{2}\left(\lambda_1+\lambda_2 +\sqrt{\lambda_1^2+\lambda_2^2-2\lambda_1\lambda_2\cos\theta}\right)}-1\right)}$$
\subsection{Reliability Measures for Other Models}
It is observed that the reliability function, $\bar F_T(\cdot)$, of the two-component series system is the same for Marshall-Olkin, Block-Basu and Sarkar's models, and is given by
$$\bar F_T(t)= e^{-t(\lambda_1+\lambda_2+\lambda_{12})}.$$
As a result, the failure rate, $r_T(t)$, and the mean residual life, $e_T(t)$, of $T$ are given by
$$r_T(t)= \lambda_1+\lambda_2+\lambda_{12}\qquad {\rm and}\qquad e_T(t)= \left(\lambda_1+\lambda_2+\lambda_{12}\right)^{-1}.$$
The reversed hazard rate function, $\mu_T(\cdot)$, is obtained as
$$\mu_T(t)=\frac{\lambda_1+\lambda_2+\lambda_{12}}{e^{t(\lambda_1+\lambda_2+\lambda_{12})}-1}.$$
In Table \ref{tab1}, we report the reliability measures for two-component series system for different models as obtained above.
	\begin{table}[h]
	\begin{center}
		\begin{tabular}{lcccl}
			\hline
			Model & $\bar{F}(t)$ & $r(t)$ & $e(t)$ & $\mu(t)$\\
			\hline
			Independent & $e^{-\lambda t}$ & $\lambda$ & $1/\lambda$ & $\frac{\lambda}{e^{\lambda t}-1}$\\&&&&\\
			Gumbel I& $e^{-(\lambda+\lambda_{12}t)t}$&$\lambda+2\lambda_{12}t$ & $e^{\delta(t)}(1-\Phi(\sqrt{2\delta(t)}))(\pi/\lambda_{12})^{1/2}$&$\frac{\lambda+2\lambda_{12}t}{e^{\lambda t+\lambda_{12}t^2}-1}$\\&&&&\\
			Gumbel II & $e^{-\lambda t}h(t)$ & $\lambda-\frac{h'(t)}{h(t)}$ &
			$\frac{\frac{1+\alpha}{\lambda}-\alpha g(t)}{h(t)}$&$\frac{\lambda h(t)-h'(t)}{e^{\lambda t}-h(t)}$\\&&&&\\
			Gumbel III&$e^{-\lambda_3 t}$&$\lambda_3$&$1/\lambda_3$&$\frac{\lambda_3}{e^{\lambda_3 t}-1}$\\&&&&\\
			Freund&$e^{-\lambda t}$&$\lambda$&$1/\lambda$&$\frac{\lambda}{e^{\lambda t}-1}$\\&&&&\\
			Marshall-Olkin& $e^{-\lambda^* t}$&$\lambda^*$&$1/\lambda^*$&$\frac{\lambda^*}{e^{\lambda ^*t}-1}$\\&&&&\\
			Block-Basu&$e^{-\lambda^* t}$&$\lambda^*$&$1/\lambda^*$&$\frac{\lambda^*}{e^{\lambda^* t}-1}$\\&&&&\\
			Cowan&$e^{-(\alpha^*/2)t}$ &$\alpha^*/2$ &$2/\alpha^*$&$\frac{\alpha^*/2}{e^{(\alpha^*/2) t}-1}$\\&&&&\\
			Sarkar&$e^{-\lambda^*t}$ & $\lambda^*$&$1/\lambda^*$&$\frac{\lambda^*}{e^{\lambda^* t}-1}$\\
			\hline
		\end{tabular}
		\caption{Reliability Measures for Two-Component Series System}\label{tab1}
	\end{center} 
\end{table}
\section{Analysis of the Errors in Reliability Measures}\label{an}
Suppose the components used to form a two-component series system actually follow some bivariate exponential distribution, and due to some reason (may be due to lack of information about the joint distribution or to get mathematical ease) we consider them to be independent. This wrong assumption incurs some errors in the inference about the system. Based on the reliability measures described in Section \ref{cal}, we analyse the behaviour of the relative errors corresponding to different bivariate exponential distributions mentioned above. 
\subsection{Error Analysis in Gumbel I Model}
In Gumbel I model the relative error in reliability function is $\left(e^{-\lambda_{12}t^2}-1\right)$, which is negative and decreasing in $t$. Note that the relative error goes from  $0$ to $-1$ when  $t$ varies from $0$ to $\infty$. It is to be mentioned  here that the horizontal line having vertical intercept $-1$ is the asymptote to the curve of relative error in reliability. 

The relative error in hazard rate is $\left(\frac{2\lambda_{12}}{\lambda}.t\right)$ which is positive and linearly increasing in $t$ having $\frac{2 \lambda _{12}}{\lambda}$ as the slope.

The relative error in mean residual life function is 
$$e^*(t)=\sqrt{\frac{\pi}{\lambda_{12}}}\;\lambda e^{\delta (t)}\left(1- \Phi\left(\sqrt{2 \delta (t)}\right)\right)-1,$$ with $\delta(t)=\lambda_{12}\left(t+\frac{\lambda}{2\lambda_{12}}\right)^2$. Since normal distribution is IFR (increasing in failure rate), Gupta and Gupta (1990) have shown, by writing $e^*(t)$ as  a function of normal failure rate, that it is decreasing. Now, it is not difficult to see that $e^*(t)$ varies from $\lambda e^{\frac{\lambda^2}{4\lambda_{12}}}\left(1-{\Phi}\left(\frac{\lambda}{\sqrt{2\lambda_{12}}}
\right)\right) \sqrt{\frac{\pi}{\lambda_{12}}}-1$ to $-1$. 

\begin{figure}
	\centering
	\includegraphics[height=8.5cm,keepaspectratio]{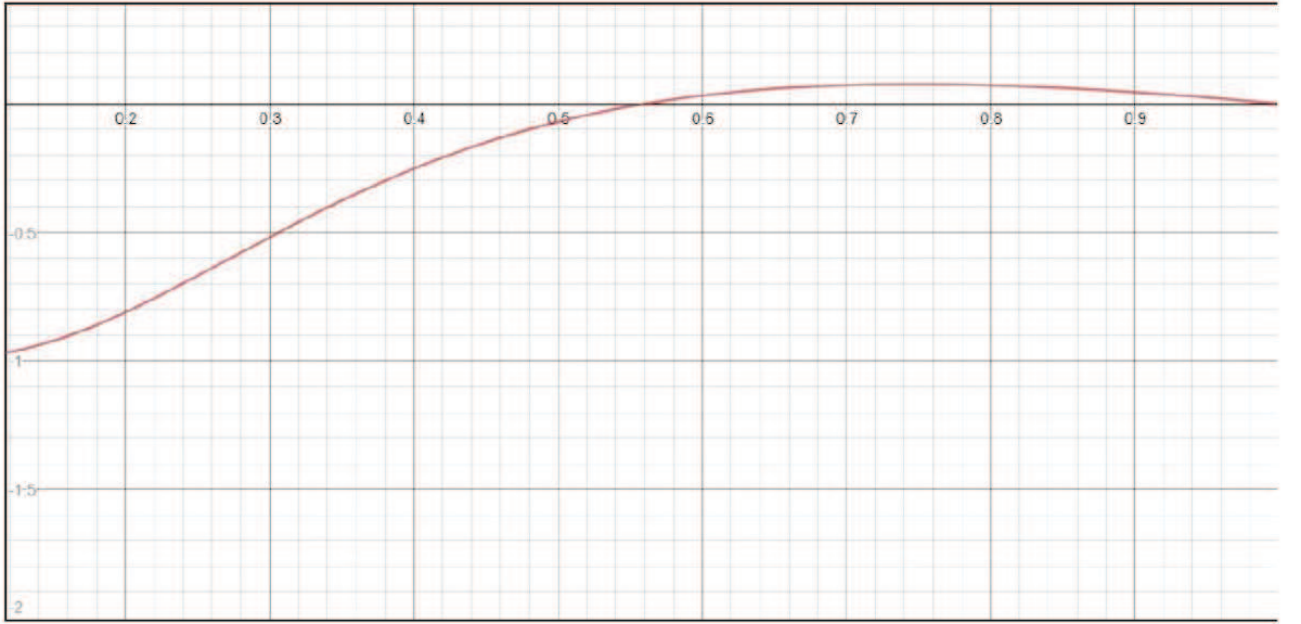}
	\caption{Relative error in RHR for Gumbel-I Model}
	\label{G1}
\end{figure}

The relative error in reversed hazard rate, given by
$$\mu^*(t)=\frac{e^{\lambda t}-1}{\lambda}\left(\frac{\lambda+2\lambda_{12}t}{e^{\lambda t +\lambda_{12}t^2}-1}\right)-1$$
with $\lambda=\lambda_1+\lambda_2$, is not monotone. To see this we take $(\lambda,\lambda_{12})=(2,1)$. Then, by writing $x=e^{-t}$, we see from Figure \ref{G1} that
$$\mu^*(t)=\frac{(1-x^2)(1-\ln x)}{x^{\ln x}-x^2}-1$$ 
is not monotone. It is observed that the function crosses the horizontal axis at the point $x=0.5616$ (approx.) which tells us that $\mu^*(t)$ crosses the horizontal axis approximately at the point $t=-\ln (0.5616)\approx 0.577$. Note that $\mu^*(t)$ increases in $t\in (0,t_0)$ and then  decreases making the horizontal line with vertical intercept $-1$ as the asymptote to the curve, where the value of $t_0$ can be obtained by solving the equation $\frac{d}{dt}\mu^*(t)=0$. The approximate value of $t_0$ (the unique value of mode) is obtained as $-\ln (0.747)\approx 0.2917$. The maximum value of $\mu^*(t)$ is $\mu^*(t_0)\approx 0.0756$.
Thus, the reliability of the Gumbel I model, under the assumption of independence, leads to over-assessment; the failure rate leads to under-assessment, whereas the mean residual life and the reversed hazard rate both give over-assessment after a certain time.

\subsection{Error Analysis in Gumbel II Model}
In Gumbel II model the relative error in reliability function is $h(t)-1$, where $h(t)$ is as given in section $3.2$. Note that $h(t)-1$ is negative and decreasing in $t$ if $-1<\alpha<0$, and positive and increasing in $t$ if $0<\alpha<1$. Note that the relative error goes from $0$ to $\alpha$ when $t$ varies from $0$ to $\infty$. It is to be mentioned here that the horizontal line having vertical intercept $\alpha$ is the asymptote to the curve of relative error in reliability.

Before we discuss the behaviour of the relative error in hazard rate in this model, we need the following lemma.
\begin{l1}\label{lem1}
	Let $$h(x)=1+\alpha (1-e^{-\lambda_1x})(1-e^{-\lambda_2x}).$$ Then, 
for $\alpha,\lambda_1, \lambda_2>0,$
$$f(x)=-\frac{h'(x)}{h(x)}\;{\rm is\; increasing\; in}\;t\in\left(\frac{2}{\lambda},\infty\right),$$ 
where $\lambda = \lambda_1 + \lambda_2$.
\end{l1} 
Proof: By differentiating $f$ with respect to $x$ we have
$$\frac{df(x)}{dx}\stackrel{sign}{=}\alpha^2 \lambda_1^2e^{-\lambda x}(e^{\lambda_2x}+e^{-\lambda_2x}-2)+\alpha^2 \lambda_2^2e^{-\lambda x}(e^{\lambda_1x}+e^{-\lambda_1x}-2)+\alpha e^{-\lambda x}(\lambda_1^2e^{\lambda_2 x}+\lambda_2^2e^{\lambda_1x}-\lambda^2),$$
which is nonnegative if $\lambda_1^2e^{\lambda_2x}+\lambda_2^2e^{\lambda_1x}-\lambda^2 \geqslant 0$, $i.e.$, if $x\geqslant \frac{2}{\lambda}$.
This gives that $f(x)$ is increasing in the region $(2/\lambda,\infty)$.\hfill$\Box$

The relative error in hazard rate is $r^*(t)=-\frac{h'(t)}{\lambda h(t)}$ which, by Lemma \ref{lem1}, is increasing in $t\geqslant \frac{2}{\lambda}$, for $0<\alpha<1$ (this condition is sufficient but not necessary) and always negative. By taking $(\alpha,\lambda_1,\lambda_2)=(0.5,1,1)$, we get that
$$r^*(t)=-\frac{h'(t)}{\lambda h(t)}=\frac{x^2-x}{2+(1-x)^2},$$
where $x=e^{-t}$. It is observed that $r^*(t)$ decreases in $t\in (0,t_0)$ where $t_0=-\ln(3-\sqrt{6})$ and then increases, making the $x$-axis an asymptote to the curve. The minimum value of $r^*(t)$ is $r^*(t_0)=-0.1124 $ (approx). The relative error in mean residual life function is 
$$e^*(t)=\frac{1+\alpha-\alpha \lambda g(t)}{h(t)}-1,$$ 
which is not monotone. To see this, we take $(\alpha,\lambda_1,\lambda_2)=(0.5,0.5,0.5)$ and see that 
$$e^*(t)=\frac{18+3x-8\sqrt x}{12+6(1-\sqrt x)^2}-1,$$
where $x=e^{-t}$. It can be observed that $e^*$ is not monotone but always positive. It is observed that $e^*(t)$ increases in $t\in(0,t_0)$, where $t_0=-\ln\left(\frac{69-9\sqrt{57}}{2}\right)$, and then decreases, making the $x$-axis an asymptote to the curve. The maximum value of $e^*(t)$ is $e^*(t_0)=0.1062 $ (approx).

\begin{figure}
	\centering
	\includegraphics[height=7cm,keepaspectratio]{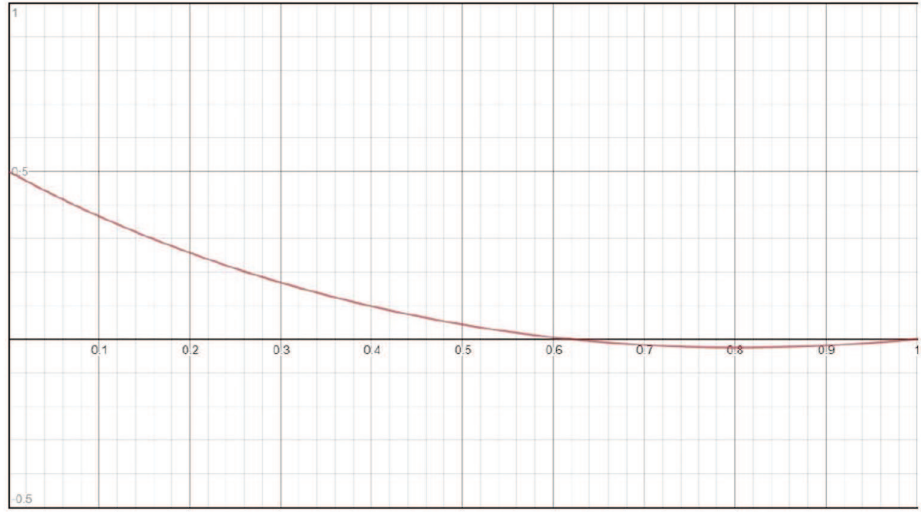}
	\caption{Relative error in RHR for Gumbel-II Model}
	\label{G2}
\end{figure}

The relative error in reversed hazard rate is given by
$$\mu^*(t)=\frac{(e^{\lambda t}-1)(\lambda h(t)-h'(t))}{\lambda (e^{\lambda t}-h(t))}-1$$ where $\lambda=\lambda_1 + \lambda_2$. It is to be noted that $\mu^*(t)$ is not monotone. To see this, we take $(\alpha,\lambda_1,\lambda_2)=(0.5,1,1)$ and get that
$$\mu^*(t)=\frac{(1-x^2)(3-3x+2x^2)}{2-2x^2-x^2(1-x)^2}-1,$$
where $x=e^{-t}$. Note from Figure \ref{G2} that $\mu^*$ is not monotone.
 It is observed that the function crosses the horizontal axis at the point $x=0.618$ (approx.) which tells us that $\mu^*(t)$ crosses the horizontal axis approximately at the point $t=-\ln (0.618)\approx 0.481$. Note that $\mu^*(t)$ decreases up to the point $t_0$ and then increases to $0.5$, making the horizontal line with vertical intercept $0.5$ as the asymptote to this curve, where the value of $t_0$ can be obtained by solving the equation $\frac{d}{dt}\mu^*(t)=0$. The approximate value of $t_0$  is obtained as $-\ln (0.8043)\approx 0.2178$. The minimum value of $\mu^*(t)$ is $\mu^*(t_0)\approx -0.0254$.

Thus, the reliability of the Gumbel II model, under the assumption of independence, leads to under-assessment or over-assessment according as $\alpha$ is positive or negative; the failure rate leads to over-assessment, whereas the mean residual life gives under-assessment and the reversed hazard rate gives under-assessment after a certain time.

\subsection{Error Analysis in Gumbel III Model}
 In Gumbel III model, by writing $\lambda_3=(\lambda_1^m+\lambda_2^m)^{\frac{1}{m}}$, we see that
 the relative error in reliability function is $e^{(\lambda - \lambda_3)t}-1$, which is positive and monotonically increases from $0$ to $\infty$ when $t$ varies from $0$ to $\infty$, since $\lambda_3 < \lambda$. The relative error in hazard rate is $\frac{\lambda_3 - \lambda}{\lambda}<0$, whereas the  relative error in mean residual life function is $\frac{\lambda - \lambda_3}{\lambda_3}>0$..
 
 Before we discuss the relative error in reversed hazard rate of this model we need the following lemma.
 \begin{l1} \label{a1}
 	For $x>0$,
 	$$f(x)=\frac{\gamma}{\beta}\left(\frac{e^{\beta x}-1}{e^{\gamma x}-1}\right)-1$$ is increasing (resp. decreasing) in $x$, provided $\beta>\gamma$ (resp. $\beta < \gamma$).
 \end{l1} 
 Proof: Note that
 $$ \frac{df(x)}{dx}=\frac{\gamma^2 e^{(\beta + \gamma)x}}{(e^{\gamma x}-1)^2}\left(\frac{e^{\gamma x}-1}{\gamma e^{\gamma x}}- \frac{e^{\beta x}-1}{\beta e^{\beta x}}\right).$$
 Let $x_0\;(>0)$ be any real number. Then
 $$f'(x_0) =\left(\frac{\gamma}{e^{\gamma x_0}-1}\right)^2 e^{(\beta + \gamma)x_0} \left( g(\gamma) - g(\beta)\right),$$ 
 where $g(t)= \frac{e^{x_0 t}-1}{te^{x_0t}}$. Taking derivative of $g(t)$ we respect to $t$, we get that 
 $$g'(t) = \frac{1+tx_0-e^{x_0t}}{t^2e^{x_0t}} < 0$$ for all $t > 0$ so that $g(t)$ is decreasing in $t$.
 Since $x_0$ is any arbitrary positive real number, the function $f$ is increasing (resp. decreasing) in $x$ when $\beta>\;({\rm resp.\;<})\; \gamma$.\hfill$\Box$
 
 The relative error in reversed hazard rate function is 
 $$\frac{\lambda_3}{\lambda}\left(\frac{e^{\lambda t}-1}{e^{\lambda_3 t}-1}\right)-1,$$ 
 which is increasing in $t$, since $\lambda > \lambda_3$ (on using Lemma \ref{a1}) and it increases from $0$ to $\infty$ when $t$ varies from $0$ to $\infty$.
 
 Thus, the reliability of this model, under the assumption of independence leads to under-assessment; the hazard rate leads to over-assessment, whereas the mean residual life and reversed hazard rate lead to under-assessment.
 
\subsection{Error Analysis in Freund's Model}
In Freund's model the relative errors in reliability function, hazard rate, mean residual life and reversed hazard rate are all zero. Thus, as far as reliability, failure rate, mean residual life or reversed hazard rate functions are concerned, the use of Freund's model and that of independent components are equivalent.
\subsection{Error Analysis in Cowan's Model}
The relative error in reliability in Cowan's model is 
\begin{equation}\label{eq1}
e^{\left(\lambda - \frac{\lambda_1+\lambda_2+\sqrt{\lambda_1^2+\lambda_2^2-2\lambda_1 \lambda_2 \cos \theta}}{2}\right)t}-1=e^{\left(\lambda - \frac{\alpha^*}{2}\right)t}-1,
\end{equation}
where $\alpha^*=\lambda_1+\lambda_2+\sqrt{\lambda_1^2+\lambda_2^2-2\lambda_1 \lambda_2 \cos \theta}$. Clearly, the expression in (\ref{eq1}) is increasing in $t$, since $\alpha ^*/2<\lambda$, and it  increases from $0$ to $\infty$ as $t$ increases, keeping $\alpha$ and $\lambda$ fixed. The relative error in the hazard rate for this model is $\frac{\alpha ^*}{2 \lambda}-1<0$, whereas the relative error in mean residual life function is $\frac{2 \lambda}{\alpha^*}-1>0$. The relative error in reversed hazard rate is 
$$\frac{\alpha^*}{2 \lambda} \left(\frac{e^{\lambda t}-1}{e^{\alpha^* t/2}-1}\right) -1,$$ 
which, by Lemma \ref{a1} and the fact that $\alpha ^*/2<\lambda$, is increasing in $t$, and it increases from $0$ to $\infty$ when $t$ varies from $0$ to $\infty$.

Thus, if the lifetimes of the two components forming a series system actually follow Cowan's bivariate exponential model, the assumption of independence leads to under-assessment of reliability, mean residual life and reversed hazard rate, whereas wrong assumption of independence leads to an over-assessment in case of failure rate.
\subsection{Error Analysis in Other Models}
In each of Marshall-Olkin, Block-Basu and Sarkar's models, the relative error in reliability function is $\left(e^{- \lambda_{12} t}-1\right)$ which monotonically decreases from $0$ to $-1$ when $t$ varies from $0$ to $\infty$. It is to be mentioned here that the horizontal line having the vertical intercept $-1$ is the asymptote to the curve of relative error in reliability. For these models, the relative error in failure rate function is $\frac{\lambda_{12}}{\lambda}>0$, whereas the relative error in mean residual life function is $- \frac{\lambda_{12}}{\lambda^*}<0$. The relative error in reversed hazard rate is 
$$\frac{\lambda^*}{\lambda} \left(\frac{e^{\lambda t}-1}{e^{\lambda^* t}-1}\right)-1,$$ 
which, by Lemma \ref{a1} and the fact that $\lambda^* > \lambda$, is decreasing in $t$, and it decreases from $0$ to $-1$, making the horizontal line, with vertical intercept $-1$, as the asymptote to this curve.

Thus, in these models, under the assumption of independence, the reliability, the mean residual life and the reversed hazard rate give over-assessment whereas the failure rate leads to under-assessment.

Table~\ref{tab2} shows the relative errors in different reliability measures for the two-component series system under the assumption of independence. 
\section{Conclusions}
Due to mathematical simplicity or otherwise, if the component lives of a two-component series system are taken to be independent when they actually follow some kind of bivariate exponential model, we may encounter over-assessment or under-assessment of the relative errors in different reliability measures. It is to be mentioned here that the relative error in reliability, hazard rate, mean residual life and reversed hazard rate is zero under the assumption of independent component lives when they actually follow Freund's bivariate exponential model, whereas it is not zero if the underlying model is any one of the models due to Gumbel (I, II, III), Marshall-Olkin, Block-Basu, Cowan or Sarkar. The analysis of the relative errors as obtained above is summerized in Table \ref{table3}, where OA stands for over-assessment and UA stands for under-assessment.
\begin{table}
	\begin{center}
		\begin{tabular}{lcccc}
			\hline
			Model & Reliability & Hazard rate & Mean residual life& Reversed hazard rate\\
			\hline
			Gumbel I &OA&UA&UA {\rm if} $t<t_0$&UA {\rm if} $t<t_1$\\&&&OA {\rm if} $t>t_0$&OA {\rm if} $t>t_1$\\&&&&\\
			Gumbel II&UA if $\alpha>0$&OA&UA&OA {\rm if} $t<t_2$\\
			              &OA if $\alpha<0$&&&UA {\rm if} $t>t_2$\\&&&&\\
			Gumbel III&UA&OA&UA&UA\\&&&&\\
			Cowan &UA&OA&UA&UA\\&&&&\\
			Marshall-Olkin&OA&UA&OA&OA\\&&&&\\
			Block-Basu&OA&UA&OA&OA\\&&&&\\
			Sarkar&OA&UA&OA&OA\\\hline 
		\end{tabular}
		\caption{Error Analysis in Two-Component Series System}\label{table3}
	\end{center}
\end{table} 
Here $t_0$ is the solution of the equation
$$\sqrt{\frac{\pi}{\lambda_{12}}}\;\lambda e^{\delta (t)}\left(1- \Phi\left(\sqrt{2 \delta (t)}\right)\right)=1$$ with $\delta(t)=\lambda_{12}\left(t+\frac{\lambda}{2\lambda_{12}}\right)^2$, $t_1$ is the solution of the equation
$$\frac{e^{\lambda t}-1}{\lambda}\left(\frac{\lambda+2\lambda_{12}t}{e^{\lambda t +\lambda_{12}t^2}-1}\right)=1$$
and $t_2$ is the solution of the equation 
$$\frac{(e^{\lambda t}-1)(\lambda h(t)-h'(t))}{\lambda (e^{\lambda t}-h(t))}=1$$
where $h(t)$ is as given in Lemma \ref{lem1}.
\begin{landscape}		
	\begin{table}
		\begin{center} 
			\begin{tabular}{lclcc}
				\hline
				Model& Reliability & Hazard Rate & Mean Residual&RHR \\
				\hline
				Independent&$0$&$0$&$0$&$0$\\&&&&\\
				Gumbel I & $e^{-\lambda_{12}t^2}-1$ & $2\lambda_{12}t/\lambda$ & $(1-\Phi(\sqrt{2\delta(t)}))(\pi/\lambda_{12})^{1/2}\lambda e^{\delta(t)}-1$&$\frac{e^{\lambda t}-1}{\lambda}\left(\frac{\lambda+2\lambda_{12}t}{e^{\lambda t+\lambda_{12}t^2}-1}\right)-1$\\&&&&\\
				Gumbel II & $h(t)-1$ & $-\frac{h'(t)}{\lambda h(t)}$& $\frac{1+\alpha-\alpha\lambda g(t)}{h(t)}-1$&$\frac{\left(e^{\lambda t}-1\right)\left(\lambda h(t)-h'(t)\right)}{\lambda\left(e^{\lambda t}-h(t)\right)}-1$\\&&&&\\
				Gumbel III& $e^{- (\lambda_3 - \lambda ) t}-1$&$(\lambda_{3}- \lambda)/\lambda$&$(\lambda -\lambda_{3})/\lambda_3 $&$\frac{\lambda_3}{\lambda}\left(\frac{e^{\lambda t}-1}{e^{\lambda_3 t}-1}\right)-1$\\&&&&\\
				Freund&$0$&$0$&$0$&$0$\\&&&&\\
				Marshall-Olkin& $e^{-\lambda_{12}t}-1$&$\lambda_{12}/\lambda$&$-\lambda_{12}/\lambda^*$&$\frac{\lambda^*}{\lambda}\left(\frac{e^{\lambda t}-1}{e^{\lambda^* t}-1}\right)-1$\\&&&&\\
				Block-Basu&$e^{-\lambda_{12}t}-1$&$\lambda_{12}/\lambda$&$-\lambda_{12}/\lambda^*$&$\frac{\lambda^*}{\lambda}\left(\frac{e^{\lambda t}-1}{e^{\lambda^* t}-1}\right)-1$\\&&&&\\
				Cowan&$e^{-(\alpha^*/2-\lambda)t}-1$&$\frac{\alpha^*}{2\lambda}-1$&$\frac{2\lambda}{\alpha^*}-1$&$\frac{\alpha^*}{2\lambda} \left(\frac{e^{\lambda t}-1}{e^{\alpha^* t/2}-1}\right)-1$\\&&&&\\
				Sarkar&$e^{-\lambda_{12}t}-1$&$\lambda_{12}/\lambda$&$-\lambda_{12}/\lambda^*$&$\frac{\lambda^*}{\lambda}\left(\frac{e^{\lambda t}-1}{e^{\lambda^* t}-1}\right)-1$\\
				\hline
			\end{tabular}
			\caption{Relative Errors in Reliability Measures under the Assumption of Independence} \label{tab2}
		\end{center}
	\end{table}
\end{landscape}

\end{document}